\documentclass[10pt]{amsart}
\usepackage{amsmath, amssymb, latexsym}

\newtheorem{theorem}{Theorem}[section]

\newtheorem{lemma}[theorem]{Lemma}
\newtheorem{corollary}[theorem]{Corollary}
\newtheorem{proposition}[theorem]{Proposition}

\newtheorem{definition}[theorem]{Definition}

\newcommand{\mn}{\par\medskip\noindent}

\newcommand{\cB}{\mathcal B}

\newcommand{\cU}{\mathcal U}
\newcommand{\dom}{{\rm dom}}
\newcommand{\Hom}{{\rm Hom}\,}
\newcommand{\Tor}{{\rm Tor\,}}
\newcommand{\Ext}{{\rm Ext\,}}
\newcommand{\Der}{{\rm Der}}
\newcommand{\Inn}{{\rm Inn\,}}
\renewcommand{\dim}{{\rm dim}}
\renewcommand{\ker}{{\rm ker}\,}
\newcommand{\coker}{{\rm coker}\,}
\newcommand{\Z}{{\rm Z}}
\def\Cz{\mathbb{C}}

\def\Zz{\mathbb{Z}}
\def\Nz{\mathbb{N}}
\newcommand{\deq}{\stackrel{ \rm def}{=}}     
\newcommand{\hot}{\mathbin{\overline{\otimes}}}

\title{$L^2$-cohomology for von Neumann algebras}
\author{Andreas Thom}
\address{Andreas Thom, Mathematisches Institut der Universit\"at G\"ottingen,
Bunsenstr. 3-5, D-37073 G\"ottingen, Germany}
\email{thom@uni-math.gwdg.de}
\urladdr{http://www.uni-math.gwdg.de/thom}

\begin{document}

\begin{abstract}
We study $L^2$-Betti numbers for von Neumann algebras, as defined by D. Shlyakhtenko and A. Connes in \cite{CS}. We give a definition
of $L^2$-cohomology and show how the study of the first $L^2$-Betti number can be related to the study of derivations with values
in a bi-module of affiliated operators. We show several results about the possibility of extending derivations from sub-algebras and about
uniqueness of such extensions. Along the way, we prove some results about the dimension function of modules over rings of affiliated
operators which are of independent interest.
\end{abstract}

\maketitle

\section{Introduction}

In \cite{CS}, A. Connes and D. Shlyakhtenko define $L^2$-Betti numbers for all complex tracial $*$-algebras $(A,\tau)$ which satisfy a 
positivity and a boundedness criterion:
\begin{enumerate}
\item $\tau(a^*a) \geq 0$ for all $a \in A$, and
\item $\forall b\in A \, \exists C >0\colon \tau(a^*b^*ba) \leq C \tau(a^*a)$ for all $a \in A$.
\end{enumerate}
\mn
Let $L^2(A,\tau)$ be the Hilbert space completion of the pre-Hilbert 
space $A$ with inner product $(x,y)_A = \tau(y^*x)$. The boundedness criterion ensures that $A$ acts as bounded operators on $L^2(A,\tau)$ and so 
the enveloping von Neumann algebra $M = W^*(A)\subset \cB(L^2(A,\tau))$ exists. The definition of the $k$-th $L^2$-Betti number of $(A,\tau)$ is now as follows:
\[\beta^{(2)}_k(A,\tau) \deq \dim_{M \hot M^o} \Tor^{A \otimes A^o}_k(M \hot M^o,A) \in [0,\infty].\]

Here, $M\hot M^o$ is naturally seen as right $A \otimes A^o$ module via the inclusion $A \otimes A^o \subset M \hot M^o$. Note that $M \hot M^o$ carries also
a left module structure over $M \hot M^o$. Since the actions commute and $\Tor_k^{A \otimes A^o}(?,A)$ is functorial, $\Tor_k^{A \otimes A^o}(M \hot M^o,A)$ inherits a left-module structure over $M \hot M^o$ with respect to which one can take a dimension. The dimension function to be used is the generalized 
dimension function of W. L\"uck, see \cite{L} on pp. 237.
\mn
The definition is modelled to give a generalization of $L^2$-Betti numbers of discrete groups. Indeed Proposition $2.3$ in \cite{CS} shows that for a
discrete group $\Gamma$ 
\[ \beta^{(2)}_k(\Gamma) = \beta^{(2)}_k(\Cz \Gamma,\tau), \quad \forall k \geq 0,\]
where $\beta^{(2)}_k(\Gamma)$ denotes the $k$-th $L^2$-Betti number of the discrete group $\Gamma$, as studied by J. Cheeger and M. Gromov in \cite{CG}.
\mn
In this article we show the vanishing of all $L^2$-Betti numbers for von Neumann algebras with diffuse center. This generalizes a result in \cite{CS}. 
We study ring-theoretic properties of $\cU(M)$ and use them to give a reasonable definition of $L^2$-cohomology for tracial algebras. After setting up the stage for
$L^2$-cohomology we concentrate on the first $L^2$-cohomology group which can be described by derivations with values in a module of affiliated operators.
The usual analytic tools to study these 'unbounded' derivations break down and we start to develop some of their basic properties.
\mn
Apart from the Introduction this article contains $5$ more sections. In Section $2$ we show the vanishing of $L^2$-Betti numbers in the case of
von Neumann algebras with diffuse center. In Section $3$ we recall and develop some algebraic properties of the ring of operators affiliated with
a finite von Neumann algebra. Here, we define $L^2$-cohomology of a tracial algebra $A$ as the Hochschild cohomology with coefficients
in a bi-module of affiliated operators. Section $4$ and $5$ are devoted to a study of the first $L^2$-cohomology group. The article ends with
concluding remarks in Section $6$.
\mn
We assume that \textit{all} von Neumann algebras appearing in this article carry a distinguished normal, faithful tracial
state. \textit{All} homomorphisms between von Neumann algebras are supposed to be unital and trace-preserving, unless otherwise stated.
All sub-algebras are unital with respect to the same unit.

Given a von Neumann algebra $M$, we denote by $M^o$ the opposite algebra. 
The notation $\cU(M)$ is used for the ring of affiliated operators, see \cite{Takesaki2}, pp. 167.
We denote by $\dim_M$ the dimension function, as defined by W. L\"uck in \cite{L}. It assigns 
to {\it any} module over $M$ an element in $[0, \infty]$ and satisfies various properties, see \cite{L}, pp. 237. We also
consider $\dim_{\cU(M)}$, the analogous dimension function for modules over the ring of affiliated operators, see \cite{Rei}.

We freely identify bi-modules over $M$ with $M \otimes M^o$-modules, left $M$-modules with right $M^o$-modules. Given a left
$M$-module $L$, put $L^o$ for the corresponding right $M^o$-module.

\section{Diffuse center}

Apart from finite dimensional algebras, the following result is so far the only 
complete computation of $L^2$-Betti numbers for von Neumann algebras. We will give a generalization 
in Theorem \ref{diffuse}.

\begin{theorem}[Connes-Shlyakhtenko, Corollary $5.4$ in \cite{CS}]\label{abel}
Let $M$ be a commutative von Neumann algebra. Then $$\beta_k^{(2)}(M) =0, \quad \forall k\geq 1.$$
\end{theorem}

It was noted in \cite{CS} that the difficult part of the preceding theorem is
the case where $M$ is diffuse, i.e. $M$ is isomorphic to $L^{\infty}(X)$ for some measure space $(X,\mu)$ without atoms.
Moreover, in Corollary $3.5$ of \cite{CS} it was shown that $\beta^{(2)}_1(M,\tau)$ vanishes in case of diffuse center.
The following theorem gives a generalization of both results.

\begin{theorem}\label{diffuse}
Let $M$ be a commutative von Neumann algebra with diffuse center. Then
$$\beta^{(2)}_k(M) = 0, \quad \forall k \geq 0.$$
\end{theorem}
\begin{proof}
We identify the centre of $M$ with $L^{\infty}[0,1]$ and consider the sub-algebra $A \subset L^{\infty}[0,1]$ 
generated algebraically by the Rademacher functions, i.e. the algebra of step-functions associated to the dyadic subdivision of $[0,1]$. The algebra
$A$ is filtered by the degree of subdivision. Denote by $A_n$ the algebra of step-functions, associated to the subdivision
$\{0,2^{-n},\dots,1\}$. 

Note that $M\otimes_{A_n}M$ is a finitely generated projective $M\otimes M^o$-module. 
Indeed, $$M \otimes_{A_n} M \cong  M \otimes M^o \left( \sum_{i=0}^{2^n-1} \chi_{[i2^{-n},(i+1)2^{-n}]} \otimes \chi_{[i2^{-n},(i+1)2^{-n}]}^o \right).$$

Furthermore, $M \otimes_A M$ is the co-limit of the $M  \otimes_{A_n} M$ in the category of $M \otimes M^o$-modules. We conclude that it is 
a flat $M \otimes M^o$-module, being a colimit of projective modules.

We consider the relative bar resolution $C_* \to M$ of $M$ by $M \otimes M^o$-modules, where $C_n = M^{\otimes_A^{n+1}}$.
A similar argument as above shows that each $C_n$ is a flat $M \otimes M^o$-module so that $C_*$ is in fact a flat resolution. 
Note that any flat resolution can be used to compute the Tor-groups and their dimensions.

First of all, we compute $$\dim_{M \hot M^o} \left( M \hot M^o \otimes_{M \otimes M^o}(M \otimes_A M)  \right)= 0.$$ Indeed, this follows
from the fact that $M \hot M^o \otimes_{M \otimes M^o} (M \otimes_{A_n} M)$ surjects onto $M \hot M^o \otimes_{M \otimes M^o}(M \otimes_A M)$ for
every $n$, and 
$$\dim_{M \hot M^o} \left( M \hot M^o \otimes_{M \otimes M^o}(M \otimes_{A_n} M)  \right)= 2^{-n}.$$

Note that $C_n$ is central with respect to $A$, so that $C_n$ is a quotient of a sum of modules of the form $M \otimes_A M$. 
By cofinality of the dimension function of L\"uck, we can conclude 
that $$\dim_{M \hot M^o}  M \hot M^o \otimes_{M \otimes M^o} C_n =0,$$ for each $n \geq 0$. This implies the claim, since $\Tor_n^{M \otimes M^o}(M\hot M^o, M)$ is a isomorphic to a sub-quotient of $M \hot M^o \otimes_{M \otimes M^o} C_n$ and hence also zero-dimensional.
\end{proof}

\section{Algebraic properties of $\cU(M)$}

The approach of Shlyakhtenko and Connes in \cite{CS} analyzed the homological side of $L^2$-invariants for tracial algebras. 
Unfortunately, so far the theory lacks a dual picture which is well understood from the point of homological algebra. In order to 
achieve this goal we introduce a new player, the algebra of affiliated operators.
\mn
In this section, let $M$ be a finite von Neumann algebra and $\cU(M)$ its ring of affiliated operators. We denote the ring
of operators affiliated with $M \hot M^o$ by $\cU^e$.

First note that for every $*$-sub-algebra $A \subset M$,
\begin{eqnarray*}
\beta_k^{(2)}(A,\tau) 
&=& \dim_{M \hot M^o} \, \Tor_k^{A \otimes A^o}(M \hot M^o,A) \\
&=& \dim_{\cU^e} \, \cU^e \otimes_{M\hot M^o} \Tor_k^{A \otimes A^o}(M \hot M^o,A)\\
&=& \dim_{\cU^e} \, \Tor_k^{A \otimes A^o}(\cU^e,A )
\end{eqnarray*}
since $\cU^e$ is flat over $M\hot M^o$ (see Theorem $8.22$ and $8.29$ in \cite{L}).
\mn
In the sequel we will show that the $\Ext$-theory for the ring of affiliated operators is much closer related to the $L^2$-homology
than this could possibly be true for the von Neumann algebra itself.

Let us first recall some algebraic properties of the ring of operators affiliated with a finite von Neumann algebra.

\begin{theorem} \label{generalities}
$\cU(M)$ is a self-injective and von Neumann regular ring.
\end{theorem}

These are classical facts and proofs can be found in \cite{Ber,Rei}. We recall that there are several characterizations of von Neumann regular rings. First
of all, they can be characterized by the algebraic property that for every element $a$ in the ring, there exists an element $b$ in the ring, such that $aba=a$ and $bab =b$.

Since we are in a $*$-algebra, we could require that the idempotents $ab$ and $ba$ are in fact projections. A computation shows that this forces
uniqueness for $b$. An element $b$ which satisfies all the mentioned properties is called Moore-Penrose 
inverse or just \textit{partial inverse} of $a$. Using the polar decomposition, one can show that the $*$-algebra
$\cU(M)$ contains partial inverses to each element. In particular, it is a von Neumann regular ring. 

A second characterization of von Neumann regular rings is in terms of homological algebra. A ring is von Neumann regular if and only if every module over it
is flat. We will use both of the characterizations. A general reference for the theory of von Neumann regular rings is \cite{Go}.
\begin{definition}
Let $A$ be a ring with involution. If $A$ contains a partial inverse to each of its elements, then $A$ is called $*$-regular.
\end{definition}
The preceding remark shows that $\cU(M)$ is a $*$-regular ring. 
\mn
A left module $L$ over a ring $R$ is called injective, if the functor $\hom_R(?,L)$ from left $R$-modules to abelian groups is exact. A ring is called 
left self-injective, if it is injective as a left module over itself. The notion of right self-injectivity is defined similarily. Since $\cU(M)$ carries an 
involution, left self-injectivity and right self-injectivity are equivalent. In the following considerations we will
exploit some of the consequences of self-injectivity of $\cU(M)$.
\mn
We denote the dual of a left $\cU(M)$-module $L$ by $$L'= \hom_{\cU(M)}(L,\cU(M)).$$ Note that $L'$ is naturally a right-module over $\cU(M)$.
Note that self-injectivity of $\cU(M)$ is equivalent to the requirement that contravariant functor that assign $L'$ to $L$ is exact from the category of left $\cU(M)$-modules
to the category of right $\cU(M)$-modules.

\begin{corollary}\label{dual}
Let $K$ be a left $\cU(M)$-module. The module $K$ is zero-dimensional if and only if $K'$ vanishes.
\end{corollary}
\begin{proof}
Let us assume that $K'$ is not zero. There exists a non-trivial homomorphism $\phi: K \to \cU(M)$, thus there exists an $x\in K$ for which
$\phi(x) \neq 0$. Let $\tilde{K}$ be the sub-module of $K$ which is generated by $x$. The $\cU(M)$-module $\phi(\tilde{K})$ is non-trivial and
finitely generated. Hence $\phi(\tilde{K})$ is projective and has positive dimension. This implies that $\tilde{K}$ and hence $K$ 
have positive dimension too.

Conversely, if $K$ has positive dimension, there exists a injection $\iota: \cU(M) p \to K$. By self-injectivity of $\cU(M)$, its dual 
$\iota': K' \to (\cU(M) p)'\neq 0$ is surjective, and hence non-trivial. Thus, $K'$ has to be non-trivial.
\end{proof}

\begin{corollary} \label{dimeq}
Let $K$ be a $\cU(M)$-module. We have an equality $$\dim_{\cU(M)} K = \dim_{\cU(M)^o} K'.$$
\end{corollary}
\begin{proof} Let us assume that $\dim_{\cU(M)} K$ is finite. Then there exists an injective map of 
$\cU(M)$-modules $\phi: \oplus_{n \in \Nz} \cU(M) p_n \to K$, such that $\sum_{n\in Nz} \tau(p_n) = \dim_{\cU(M)} K$. By self-injectivity
of $\cU(M)$ the $\phi': K' \to \prod {\cU(M) p_n}'$ is surjective. Furthermore, $\ker{\phi'} = (\coker{\phi})'$ is of dimension zero, since
$\coker \phi$ is of dimension zero and Corollary \ref{dual}.
Thus, the claim follows from following chain of the equalities:
$$\dim_{\cU} K' = \dim_{\cU(M)^0} \prod_{n \in \Nz} \cU(M)^o p_n = \sum_{n \in \Nz} \tau(p_n) = \dim_{\cU(M)} K.$$
For the equality in the middle, we used Theorem $6.18$ in \cite{L}.
\end{proof}

Note that the compatibility under duality implies a nice behaviour of the dimension function 
with respect to arbitrary limits in the category of $\cU(M)$-modules. This is far from
being true for $M$-modules. Compare with Theorem $6.18$ in \cite{L}.

\begin{theorem}\label{duality}
Let $A\subset \cU(M)$ be a $*$-sub-algebra and let $K$ be any $A$-bi-module. There exists a natural map
$$\Ext^k_{A\otimes A^o}(K,\cU^e) \to \Tor_k^{A \otimes A^o}(\cU^e,K)'$$
which is an isomorphism of right $\cU^e$-modules.
\end{theorem}
\begin{proof}
This follows from the fact (see e.g. Proposition $2.6.3$ in \cite{W}), that for any $A \otimes A^o$-module $L$, there is a natural isomorphism of
right $\cU^e$-modules
$$\Hom_{A\otimes A^o}(L,\cU^e) \to \Hom_{\cU^e}(\cU^e \otimes_{A \otimes A^o}  L , \cU^e),$$
and the fact that duality is exact, i.e. $\cU^e$ self-injective. 

Indeed, let $C_*$ is a projective resolution of $K$. The homology
of $ \cU^e \otimes_{A \otimes A^o} C_*$ computes the Tor-groups. Dualizing is exact, so that taking homology first and then dualizing
is the same as dualizing first and taking then cohomology. This implies that the duals of the Tor-groups can be computed as the cohomology
of the complex $\Hom_{\cU^e}(\cU^e \otimes_{A \otimes A^o} C_* ,\cU^e)$. By the observation above, this complex of right $\cU^e$-modules
is naturally isomorphic to $\Hom_{A \otimes A^o}(C_*,\cU^e)$. The cohomology of this complex computes by definition 
the Ext-groups and the proof is finished.
\end{proof}

By Lemma $9.1.3$ in \cite{W}, $\Ext^k_{A\otimes A^o}(A,\cU^e)$ identifies with the $k$-th Hochschild cohomology of $A$ with coefficients
in the $A$-bi-module $\cU^e$. Similarily, as noted in \cite{CS}, $\Tor_k^{A \otimes A^o}(\cU^e,A)$ can be identified with the $k$-the Hochschild
homology of $A$ with coefficients in the $A$-bi-modules $\cU^e$. Moreover, it can be checked that the induced map
$$H^k(A, \cU^e) \to \left(H_k(A,\cU^e)\right)'$$ is the well-known pairing in Hochschild homology.

\begin{definition}
The $k$-th $L^2$-cohomology of the tracial $*$-algebra $(A,\tau)$ (satisfying the positivity and boundedness criterion mentioned in the introduction) 
is defined to be the right $\cU^e$-module
\[\Ext^k_{A\otimes A^o}(A,\cU^e) = H^k(A,\cU^e).\]
\end{definition}

\begin{corollary} \label{duality2}
Let $A$ be a $*$-sub-algebra of $\cU(M)$. The natural map
$$H^k(A,\cU^e) \to H_k(A,\cU^e)' \cong \left(H_k(A, M \hot M^o) \otimes_{M \hot M^o} \cU^e \right)'$$ is an isomorphism of right $\cU^e$-modules.
\end{corollary}

This corollary is useful for several reasons. First of all, by Corollary \ref{dimeq}, it allows to compute the $L^2$-Betti numbers of the tracial algebra
$(A,\tau)$ as the dimensions of
certain Hochschild cohomology groups which are im some cases more accessible than the Hochschild homology groups. Secondly, a vanishing result of
$L^2$-Betti numbers can be translated in a real vanishing result of cohomology groups by Corollary \ref{dual}. This is of course
very useful since the study of zero-dimensional modules is circumvented.

Of course, it is natural to consider topological versions of Hochschild cohomology since the ring and the coefficient module carry natural topologies. 
In the sequel we will do this only in the case of the first $L^2$-cohomology group.

\section{First $L^2$-cohomology for von Neumann algebras}

\subsection{Derivations}

Let $A$ be a unital $\Cz$-algebra and $L$ be a $A$-bimodule. Denote by $\Der(A,L)$ the set of $\Cz$-linear derivations with values in the bi-module $L$, i.e.
$\Cz$-linear maps $\delta\colon A \to L$, which satisfy $$\delta(ab) = a\delta(b) + \delta(a)b,$$ for all $a,b \in A$. A derivation $\delta \in \Der(A,L)$ 
is called inner if there exists $\xi \in L$, such that $$\delta(a) = [a,\xi], \quad \forall a \in A.$$ We use the notation $\Inn(A,L)$ for the set of inner
derivations.
\mn
By Remark $1.5.2$ in \cite{Lo}, the first Hochschild cohomology of a ring $A \subset \cU$ with coefficients in the bi-module $\cU^e$ is described
as $$H^1(A,\cU^e) = \frac{\Der(A,\cU^e)}{\Inn(A,\cU^e)}.$$ Here, we use the left $A \otimes A^o$ module structure of $\cU^e$ which is 
induced by left-multiplication of $A \otimes A^o \subset \cU^e$.
Since $\cU^e$ carries a commuting right $\cU^e$-module structure, the vector spaces 
$\Der(A,\cU^e)$ and $\Inn(A,\cU^e)$ are right $\cU^e$-modules. 

Furthermore, if $A$ is an ultra-weakly dense $*$-sub-algebra of $M$, there is an exact sequence of right $\cU^e$-modules as follows:
\[ 0 \to {\Z_A}(\cU^e)\to  {\cU^e} \to {\Der}(A,\cU^e)\to H^1(A,\cU^e) \to 0. \]
Here, $\Z_A(\cU^e)$ denotes the centre of $\cU^e$ with respect to $A$, i.e. 
\[\Z_A(\cU^e) = \{ \xi \in \cU^e\colon (a \otimes 1)\xi = (1 \otimes a^o) \xi, \forall a \in A  \}.\]

It was shown in \cite{CS} that $$\dim_{M \hot M^o} \Z_A(M \hot M^o) = \beta_0^{(2)}(M).$$ The corresponding result
in our case (with an identical proof) says that $$\dim_{\cU^e} \Z_A(\cU^e) = \beta_0^{(2)}(M).$$
The exactness  of the above sequence and the additivity of the dimension function imply the following proposition.
\begin{proposition} Let $A$ be an ultra-weakly dense $*$-sub-algebra of $M$.
$$\dim_{\cU^e} \Der(A,\cU^e)= \beta_1^{(2)}(A)- \beta_0^{(2)}(M) +1.$$
\end{proposition}

\subsection{Group von Neumann algebras}

It is hard to distinguish von Neumann algebras. One reason is, that in many cases, it turned out that completely different constructions lead
to only one von Neumann algebra. Let us mention only that all the group von Neumann algebras of i.c.c. amenable groups are isomorphic.

Although already von Neumann himself showed that the group von Neumann algebras of free groups are not 
isomorphic to ones coming from amenable groups, it remained open, whether non-isomorphic free groups lead to non-isomorphic group von 
Neumann algebras. It is well-known that
$\Cz F_n \not\cong \Cz F_m$ unless $n = m$. (E.g. $\hom(\Cz F_n,\Cz) = (\Cz^\times)^n$ and $\pi_1((\Cz^\times)^n)= \Zz^n$.) 
There are many ways to distinguish these group algebras and one is the first $L^2$-Betti-number:
\[ \beta^{(2)}_1(\Cz F_n,\tau) = n-1.\]

The authors of \cite{CS} formulate the following equality not as a conjecture, but the philosophy the lead to the definitions in
\cite{CS} certainly implies that some version of the following equality should hold.
\[ \beta^{(2)}_k(\Cz \Gamma,\tau) = \beta^{(2)}_k(L\Gamma,\tau), \quad \forall k \geq 0.\]

It is clear that such a result would immediately answer the question about non-isomorphism of free group factors affirmatively.
Coming back to the results in the preceding section, we should compare the first $L^2$-cohomology of a group $\Gamma$ with the 
first $L^2$-cohomology of the group von Neumann algebra $L\Gamma$. This leads to the following two questions:

\begin{enumerate}
\item Which derivations $\Delta: \Cz \Gamma \to \cU(L (\Gamma \times \Gamma^o))$ extend to the group von Neumann algebra $L\Gamma$? 
\item Are the extensions unique?
\end{enumerate} 

We are not able to answer any of the above questions completely. However, partial results are obtained in the remaining sections. Note that
there is no automatic continuity result for derivations $L \Gamma \to \cU(L (\Gamma \times \Gamma^o))$, so that a derivation need not
be determined by its values on a dense subset. Passing to a continuous version of $L^2$-cohomology would circumvent this problem and
the questions should be reformulated accordingly.

\subsection{Topologies on $\cU(M)$}
Let us first set up some definitions. Let $M$ be a finite von Neumann algebra. On $\cU(M)$ we consider various topologies. 
First of all, there is the topology of convergence in measure, which
is the vector space topology with generating open sets 
$$V_\varepsilon = \{x \in \cU(M); \exists p=p^2=p^* \in M, \tau(p^\perp)\leq \varepsilon, \|xp\| \leq \varepsilon\}.$$
Secondly, we consider the topology which is induced by the rank metric. The rank metric is described as follows:
$$d(x,y) = N(x-y) = \dim_{\cU(M)} \, \cU(M)(x-y).$$ 
Similarily as above, its topology is generated by translates of
$$B^d_{\varepsilon}(0)= \{x \in \cU(M); \exists p=p^2=p^* \in M, \tau(p^\perp)\leq \varepsilon, xp = 0\}.$$
Clearly, the rank metric does not induce a vector space topology.

On $M$ we consider the topology, whose closed sets are the ones which are closed 
in the topology of convergence in measure, when restricted to bounded sets of $M$. We call this topology the topology of bounded convergence in measure.

\subsection{Completion with respect to the rank metric}

\begin{lemma} \label{support} Let $A$ be a $*$-regular sub-ring of $\cU(M)$. 
All derivations of $A$ with values in the bi-module $\cU^e$ are Lipschitz with constant $2$ with respect to the rank metric. \end{lemma}
\begin{proof} Let $\delta: A \to \cU^e$ be a derivation. 
Let $\varepsilon\geq 0$ be given. Let $a \in A$ with $N(a) \leq \varepsilon/2$. 
This is equivalent to $\tau(e) \leq \frac12$, where $e$ is the support projection of $a$. Since $A$ is $*$-regular, $e \in A$.

The following calculation implies the claim.
\begin{eqnarray*}
N(\delta(a)) &=& N(\delta(ae)) \\
&=&    N(\delta(a)e + a \delta(e)) \\
&\leq& N(e) + N(a) =2 N(e) \leq \varepsilon
\end{eqnarray*} \end{proof}

\begin{theorem} \label{pedersen}
Let $M$ be a von Neumann algebra and let $A\subset M$ be an ultra-weakly dense sub-$C^*$-algebra . 
Under these conditions, $A$ is dense in $M$ with respect to the rank metric.
\end{theorem}
\begin{proof} This is Theorem $2.7.3$ from \cite{Ped}, which can be thought of as being a generalization of Lusin's theorem. \end{proof}

\begin{proposition} \label{ext}
Let $M$ be a von Neumann algebra and let $A \subset M$ be a $*$-sub-algebra which is closed
under smooth functional calculus for self-adjoint elements. 
Every derivation $\Delta : A \to \cU^e$ has a unique extension to the completion of $A$ with respect to
the rank metric.
\end{proposition}
\begin{proof} Every element $a \in A$ can be written as a product $(a^*a)^{\frac14} \cdot (a^*a)^{-\frac14}a$. 
If $a$ is small in the rank metric then both factors are small in the rank metric. As in the proof of Lemma \ref{support}, 
it follows by the Leibniz rule that also $\Delta(a)$ is small in the rank metric. Furthermore, by the same computation 
$\Delta$ is certainly Lipschitz with constant $2$. This implies that $\Delta$ extends uniquely to the completion.
\end{proof}

\begin{theorem} \label{3algebras} 
Let $M$ be a finite von Neumann algebra. Let $A \subset M$ be a sub-$C^*$-algebra. Let $W^*(A)$ be its ultra-weak closure and let $\cU(W^*(A))$ its
ring of affiliated operators. Then the induced restriction 
maps $$H^1(\cU(W^*(A)),\cU^e) \longrightarrow H^1(W^*(A),\cU^e) \longrightarrow H^1(A,\cU^e)$$ are isomorphisms of $\cU^e$-modules.
\end{theorem}
\begin{proof}
Since the inclusion $\cU(W^*(A) \hot W^*(A)^o)  \subset \cU(M \hot M^o)$ is a flat ring extension, the general case follows from the case 
$W^*(A)=M$ and Corollary \ref{duality2}. 
By Theorem \ref{pedersen}, $A$ is dense in $M$ with respect to the rank metric. Also, $M$ is dense in $\cU(M)$ with respect to rank metric. The claim
follows from Proposition \ref{ext}.
\end{proof}

The following corollary is a first step towards a relation between the $L^2$-Betti number of a von Neumann algebra and the $L^2$-Betti numbers of
its sub-algebras. 

\begin{corollary}
Let $M$ be a finite von Neumann algebra and let $A \subset M$ be an ultra-weakly dense sub-$C^*$-algebra. The following equality holds: 
\begin{equation*}\beta^{(2)}_1(A,\tau) = \beta^{(2)}_1(M,\tau). \end{equation*}
\end{corollary}

\subsection{Exploiting $*$-regularity}

The following theorem shows that $*$-regularity for abelian algebras implies vanishing of the first $L^2$-cohomology group. The core of the proof is
Theorem $5.1$ in \cite{CS}.

\begin{theorem} \label{inner}
Let $M$ be a finite von Neumann algebra and let $A$ be a $*$-regular abelian sub-algebra of $\cU(M)$. All derivations $\Delta: A \to \cU^e$ are inner.
\end{theorem}
\begin{proof}
We have to compute $H^1(A,\cU^e)$. First of all, it is is the dual of $\Tor_1^{A \otimes A^o}(\cU^e,A)$ by Corollary
\ref{duality2}. Using Corollary \ref{dual}, it suffices 
to show that $$\dim_{\cU^e} \Tor_1^{A \otimes A^o}(\cU^e,A)=0.$$ 

Denote the closure of $A \subset \cU(M)$ in the topology of convergence in measure by $B$. Clearly, $\tilde{B}=B \cap M$ is an abelian finite von Neumann algebra and
$\cU(\tilde{B}) = B$. By abuse of notation, we set $B \hot B^o = \cU(\tilde{B} \hot \tilde{B})$.
Since $A$ is $*$-regular, the inclusion $A \otimes A^o \subset B \otimes B^o$ is a flat ring extension. We can conclude that
\begin{equation} \label{torsion}
\Tor_k^{B \otimes B^o}(B \hot B^o, (B \otimes B^o) \otimes_{A \otimes A^o} A ) \cong \Tor_k^{A \otimes A^o}(B \hot B^o,A), \quad \forall k \geq 0.
\end{equation}
\mn
We finish the proof by a computation. 
\begin{eqnarray*}
0 &=&\dim_{\tilde{B} \hot \tilde{B}^o} \Tor_1^{\tilde{B} \otimes \tilde{B}^o}(\tilde{B} \hot \tilde{B}^o,(B \otimes B^o) \otimes_{A \otimes A^o} A) \quad \mbox{(by Thm. $5.1$ in \cite{CS})}  \\
  &=&\dim_{B \hot B^o} \Tor_1^{B \otimes B^o}(B \hot B^o,(B \otimes B^o) \otimes_{A \otimes A^o} A)  \\
  &=& \dim_{B \hot B^o} \Tor_1^{A \otimes A^o}(B \hot B^o,A) \quad \mbox{(by Equation \ref{torsion})}\\
  &=& \dim_{\cU^e} \, \cU^e \otimes_{B\hot B^o} \Tor_1^{A \otimes A^o}(B \hot B^o, A) \\
  &=& \dim_{\cU^e} \, \Tor_1^{A \otimes A^o}(\cU^e,A)
\end{eqnarray*} \end{proof}

\begin{corollary} \label{innersub}
Let $M$ be a finite von Neumann algebra and let $A \subset M$ be an abelian ultra-weakly closed $*$-sub-algebra. Every derivation $\Delta: A \to \cU^e$
is inner.
\end{corollary}
\begin{proof}
The algebra $\cU(A)$ is a abelian $*$-regular sub-algebra of $\cU(M)$. Since $\Delta$ extends uniquely to $\cU(A)$ by Theorem \ref{3algebras}
and is inner by Theorem \ref{inner} it has to be inner on $A$.
\end{proof}

\section{Derivations}
\subsection{Uniqueness of extensions}
In this section we will address several questions concerning the uniqueness of extensions of derivations from $*$-sub-algebras. Theorem \ref{inner} which in
essence relies on Theorem $5.1$ in \cite{CS} will turn out to be quit useful in that respect.

\begin{definition}
Let $M$ be a von Neumann algebra and let $A$ be  ultra-weakly dense $*$-sub-algebra. 
Denote by $M^A$ the smallest $*$-sub-algebra of $M$, which satisfies the following properties
\begin{enumerate}
\item $A \subset M^A$, and
\item for each $n \in \Nz$, $M_n(M^A)$ is closed under taking ultra-weak closures of abelian $*$-sub-algebras.
\end{enumerate}
\end{definition}

The following theorem is the first result which relates derivations on $M^A$ to derivations on $A$. It will be complemented by
Theorem \ref{complement}.

\begin{theorem} \label{main} 
Let $M$ be a von Neumann algebra and let $A$ be an ultra-weakly dense $*$-sub-algebra of $M$ such that
$A \cap \Z(M)$ is ultra-weakly dense in $\Z(M)$. The natural restriction map 
$$(\iota_{A \subset M^A})^*\colon H^1(M^A,\cU^e) \to H^1(A,\cU^e)$$ is an injective map of $\cU^e$-modules.
\end{theorem}

\begin{proof}
Let $\Delta: M^A \to \cU^e$ be a derivation that is inner on $A$, i.e. $\Delta(a) = [a,\xi]$, for all $a \in A$. By substracting the
inner derivation, we might assume that $\Delta$ vanishes identically on $A$. We need to show that $\Delta$ vanishes on $M^A$ as well. 
Without loss of generality, we can also assume that $\Delta$ is $*$-preserving, i.e. 
$\Delta(x^*) = \Delta(x)^\#$, where $(x \otimes y^o)^\#=y^* \otimes (x^o)^*$. 
Indeed, the real part of $\Delta$, which is defined by $\Delta^r(x) = \frac12(\Delta(x) + \Delta(x^*)^\#)$, and the imaginary part, which 
is defined by $\Delta^i(x) = -\frac{i}2(\Delta(x) - \Delta(x^*)^\#)$ both vanish on $A$ and are $*$-preserving. Vanishing of $\Delta^r$ and
$\Delta^i$ on $M^A$ imply vanishing of $\Delta$, since $\Delta = \Delta^r + i \Delta^i$.

Let $\tilde{A}= \{x \in M^A; \Delta(x)=0\}$ be the subset of $M^A$, on which $\Delta$ vanishes. Clearly,
$\tilde{A}$ is a $*$-sub-algebra of $M^A$. Let us note that $\Z(M) \subset \tilde{A}$. Indeed,
since $\Delta$ is inner on $\Z(M)$ by Corollary \ref{innersub} and is zero on an ultra-weakly dense $*$-sub-algebra by assumption, 
it has to vanish on $\Z(M)$.

Every maximal commutative subalgebra of $M_n(\tilde{A})$ is ultra-weakly closed. Indeed, let $A_0$ be a maximal commutative $*$-sub-algebra of $M_n(\tilde{A})$.
The derivations $M_n(\Delta)$ vanishes on $A_0$ and is defined on its weak closure by assumption. 
By Theorem \ref{inner}, $M_n(\Delta)$ is inner on the weak closure as well and so $M_n(\Delta)$ has to vanish there, again by continuity. 
By maximality, the algebra $A_0$ is ultwa-weakly closed. This shows that $\tilde{A}=M^A$.
\end{proof}

We showed that there is a 'large' algebra $M^A$, for which extension of derivations from $A$ are unique if they exist.
It is of course important to decide whether the inclusion $M^A \subset M$ is strict. 
However, we made no progress in answering this question.

\subsection{Closable and totally discontinuous derivations}

In this section we prove that there is a dichotomy between totally discontinuous derivations and closable derivations. Let us first
give some definitions. We will not be able to conclude much from the results in this section. However, they clarify the difference between
continuous and arbitrary derivations to a certain extend.

\begin{definition}
Let $M$ be a von Neumann algebra and let $A \subset M$ be an ultra-weakly dense $*$-sub-algebra. 
A derivation $\Delta: A \to \cU^e$ is called closable, if it is closable from the norm-topology to the topology of convergence in measure.
\end{definition}

Throughout this section, $\cU^e$ is considered with the topology of convergence in measure. It is clear, that closable derivations 
$\Delta: A \to \cU^e$ are of special interest, since these are the ones which can be extended to larger domains of
definition by continuity. In the sequel we will analyze the obstruction to closability.

\begin{definition} Let $A\subset M$ be an ultra-weakly dense $*$-sub-algebra. Consider a derivation $$\Delta\colon A \to \cU^e.$$
\begin{itemize}
\item 
The range $r(\Delta)$ of $\Delta$ is defined to be the smallest left $M \hot M^o$-sub-module of $\cU^e$, which contains the image of $\Delta$.
\item
The range of discontinuity $d(\Delta)$ of $\Delta$ is defined to be the set:
\[ d(\Delta)=\{\xi \in \cU^e; \exists (x_n)_{n \in \Nz} \subset A \colon (\|x_n\| \to 0) \wedge (\Delta(x_n) \to \xi)\} 
\]
\item
The derivation is called totally discontinuous of $r(\Delta)=d(\Delta)$.
\end{itemize}
\end{definition}

Clearly, a derivation $\Delta$ is closable if and only if $d(\Delta)=0$. The following theorem separates the closable from the totally discontinuous part
of a derivation.

\begin{theorem} Let $A \subset M$ be an ultra-weakly dense $*$-sub-algebra and let $\Delta\colon A \to \cU^e$ be a derivation.
The inclusion $d(\Delta) \subset r(\Delta)$ is an inclusion of left $\cU^e$-submodules of $\cU^e$. Furthermore, 
there exists a sequence $(p_n)_{n \in \Nz}$ of mutually orthogonal projections in $M\hot M^o$ such that 
$$\Delta = \Delta (\sum_n p_n)^\perp + \sum_n \Delta p_n,$$
\[r(\Delta p_n ) = d(\Delta p_n),\, \forall n \in \Nz \quad \mbox{and} \quad d(\Delta (\sum_n p_n)^\perp) = 0.\]
In particular, $\Delta p_n$ is totally discontinuous for all $n \in \Nz$ and $\Delta (\sum_n p_n)^\perp$ is closable.
\end{theorem}
\begin{proof}Let us first show that $d(\Delta)$ is indeed a left $\cU^e$-module. We first show that it is a left module over $A \otimes A^o$. Indeed,
if $x_n \to 0$ and $\Delta(x_n) \to \xi$, then $\Delta(ax_n) = (a \otimes 1)\Delta(x_n) + (1 \otimes x_n^o)\Delta(a)\to (a \otimes 1) \xi$ and
similarily $\Delta(x_n a) \to (1 \otimes a^o)\xi$. Secondly, observe that $d(\Delta)$ is closed. Indeed, if $\xi_n \to \xi$ and 
$\xi_n \in d(\Delta)$, then there exist sequences $(x_{m,n})_{m \in \Nz}$ which converge to $0$, and $\Delta(x_{m,n}) \to \xi_n$ as $m$ tends to
$\infty$. Picking suitable sub-sequences, we may assume that both sequences converge uniformly in $n$. 
This implies that $x_{n,n}$ tends to zero in norm and that $\Delta(x_{n,n}) \to \xi$ in the topology of convergence in measure.

Being closed, $d(\Delta)$ immediatly gets a left $\cU^e$-module structure and is of the form $\cU^ep$ for some projection $p \in M\hot M^o$.
Note, that $r(\Delta p) = d(\Delta p)$. We continue this process with $\Delta p^{\perp}$ in place of $\Delta$. The trans-finite iteration of this
process ends. Since the traces of the projections we find are positive and sum up to a real number, the process actually ends after a countable number
of steps. This finishes the proof.\end{proof}

The following facts are recalled for the convenience of the reader and to motivate the next section.

\begin{theorem} \label{closcont}
Every closable derivation $\Delta: M \to \cU^e$ is continuous.
\end{theorem}
\begin{proof}
Both spaces carry complete metrizable vector space topologies. Hence a closed graph theorem holds, see e.g. Thm. $1.2.5$ in \cite{Rudin}.
\end{proof}

Let $A \subset M$ be an ultra-weakly dense $*$-sub-algebra and let $\Delta: A \to \cU^e$ be a closable derivation. Let $G(\Delta) \subset M \times \cU^e$ 
be the graph of $\Delta$. Its closure defines a linear operator $\Delta'$ defined on $$\dom(\Delta')=\pi_1\left(\overline{G(\Delta)}\right) \subset C^*(A).$$

Consider the mapping $\pi_1: \overline{G(\Delta)} \to C^*(A)$. If $\dom(\Delta')$ is of second category, it follows from Thm. $2.11$ in \cite{Rudin}
that $\dom(\Delta')=C^*(A)$ and that $\pi_1$ is open. Hence, the linear operator $\Delta'$ is defined on $C^*(A)$ and continuous.

Having this observation in mind, it is desirable to show that $\dom(\Delta')$ is rather large. We will give partial results in the next section. However, we are
unable to show that it is of second category.

\subsection{Existence of extensions}

In this section we want to analyze the chance of extending the derivation from an ultra-weakly dense $*$-sub-algebra $A \subset M$. 
We will focus on closable and $*$-preserving derivations. Recall that $*$-preserving means that $\Delta(x^*) = \Delta(x)^\#$, where
$(x \otimes y^o)^\# = y^* \otimes (x^*)^o$.
\mn
Let us first set up some notation.
For a positive operator $1 \leq f = f^\#\in \cU^e$, we consider the vectorspace
\[ \cU^e_f= \{a \in \cU^e, \|a f^{-\frac12}\| < \infty \}.\]
It is obvious that $\cU^e_f$ are left-modules over $M \hot M^o$. Furthermore, $f\leq g$ implies that
$\|a f^{-\frac12}\| \geq \|a g^{-\frac12}\|$ so that the $\cU^e_f$ form a directed system of $M \hot M^o$-modules
and \[\cU^e = \bigcup_{f \geq 1}  \cU^e_f.\]
Indeed, for each $a \in \cU^e$, we have that $a \in \cU^e_{1+a^*a + (a^*a)^\#}$.

Note that the $\cU^e_f$ are in fact Banach spaces with the norm
$$\|a\|_{(f)} = \|a f^{-\frac12}\|.$$

Let now $A$ be a $*$-sub-algebra of $M$ and let $\Delta: A \to M \hot M^o$ be a derivation. (The range of the derivation is choosen to be $M \hot M^o$ 
for sake of simplicity. 
If $A$ is countably generated as $\Cz$-algebra, then for every derivation $\Delta$ with range $\cU^e$, there exists a projection in $M \hot M^o$ with arbitrarily small
trace so that $\Delta p^{\perp}$ takes values in $M \hot M^o$.)
On $A$ we consider the norm $\|.\|_f$, which is defined as 
\[\|a\|_f = \|a\| + \|\Delta(a)\|_{(f)}.\]

Note that the Leibniz rule implies that $\|.\|_f$ is sub-multiplicative. The completion of $A$ with respect to $\|.\|_f$ is called $A_f$. Moreover, its 
directed union is denoted by
\[ A_\infty = \bigcup_{f \geq 1} A_f.\] We also consider the completion of $A_\infty$ 
with respect to the rank metric, it is called $A_{\overline{\infty}}.$ Similarily,
the completion of $A_f$ with respect to the rank metric is denoted by $A_{\overline{f}}.$ 
Since we assumed $\Delta$ to be closable, we can identify the algebras $A_{f}, A_{\overline{f}}, A_\infty$ and $A_{\overline{\infty}}$ 
with $*$-sub-algebras of $\cU(M)$.

The closure of the graph of $\Delta$ defines a canonical extension of $\Delta$ to $A_{\overline{\infty}}$ which we denote by 
$\Delta': A_{\overline{\infty}} \to \cU^e$. Moreover, $\Delta'$ is the closure of $\Delta$. Indeed, if $x_n \to x$ and
$\Delta(x_n) \to \xi$ in measure, then there exists (by the generalized) Egoroff's Theorem a 
sub-sequence $(x_{n_k})_{k \in \Nz}$ so that $\Delta(x_{n_k}) \to \xi$ almost
surely. In particular, there exists some invertible positive $f \in \cU^e$ as above, such that $\Delta(x_{n_k})f \to \xi f$ in norm.
\mn
For each $f$, $\|.\|_f$ is a differential norm on $A$ in the sence of \cite{CuntzBlackadar}. Hence $A_f \subset M$ is closed under smooth functional calculus for
normal elements and holomorphic calculus for arbitrary elements. In particular, the spectral radius in $A_f$ is the 
same as in $M$. We want to show that $A_{\overline{\infty}}$ is a $*$-regular sub-algebra of $\cU(M)$.

\begin{lemma}
For each element $a \in A_{\overline{\infty}}$, the partial inverse of $a$ is also contained in $A_{\overline{\infty}}$.
\end{lemma}
\begin{proof}
Assigning to each element its partial inverse is a self-map of $\cU(M)$ which is continuous with respect to the rank metric. Therefore, it
suffices to show that partial inverses to elements $a \in A_{\infty}$ exist in $A_{\overline{\infty}}$. Furthermore, it is sufficient
to prove the result for positive operators. Indeed, if $b$ is a partial inverse of $a^*a$, then $a^* b$ will be a partial inverse to $a$.

Let $a$ be a positive operator in $A_f$, for some $f \in \cU^e$ as above. Let $g_n \in C^{\infty}(-1,\infty)$ be a sequence of
smooth functions such that $g_n$ is bounded by $2n$, $g_n(0) =0$ and $g_n(t) = t^{-1}$ for all $t \geq n^{-1}$. Since $A_f$ is closed under smooth
functional calculus for self-adjoint elements, $g_n(a)$ exists in $A_f$. Furthermore, $g_n(a)$ converges with respect to the rank metric 
to the partial inverse of $a$, as $n$ tends to infinity. This shows that a partial inverse to $a$ exists in $A_{\overline{\infty}}$.
\end{proof}
Let $B \subset A_{\overline{\infty}}$ be a maximal abelian $*$-sub-algebra. It follows that $B$ is $*$-regular and by Theorem \ref{inner}
we conclude that $\Delta'$ is inner on $B$, i.e.
\[\Delta'(b) = (b \otimes 1 - 1 \otimes b^o) \xi\]
for some $\xi \in \cU^e$ and for all $b \in B$. By continuity, the closure of $B$ in the topology of bounded convergence in measure 
is contained in  $A_{\overline{\infty}}$.

In particular, we have shown that the algebra $A_{\overline{\infty}}$ is closed under measurable functional calculus for normal elements.

\begin{theorem} \label{complement}
Let $M$ be a finite von Neumann algebras and let $A$ be ultra-weakly dense $*$-sub-algebra of $M$. 
Let $\Delta: A \to \cU^e$ be a closable derivation. There is a unique extension of the derivation $\Delta$ to the algebra $M^A$. 
\end{theorem}
\begin{proof}
This follows immediately from Theorem \ref{main} and the observation that $M^A \subset A_{\overline{\infty}}$. Indeed,
let $B \subset M_n(A_{\overline{\infty}})$ be a maximal abelian sub-algebra. We have to show that the ultra-weak closure
is also contained in $M_n(A_{\overline{\infty}})$. The assertion has been proved for the particular case $n=1$. 

Note that $M_n(A)_{\infty} = M_n(A_{{\infty}})$ since every positive operator in $M_n(\cU^e)$ is dominated by an operator of
the form $1_{M_n \Cz} \otimes a$. Moreover, a matrix converges in rank metric if and only each of its entries converges in the
rank metric. This shows that 
$$M_n(A)_{\overline{\infty}} = M_n(A_{\overline{\infty}}).$$

Now, since $M_n(A_{\overline{\infty}})$ is $*$-regular it follows that $B$ is $*$-regular and hence that 
$M_n(\Delta)$ is inner on $B$. As before, we conclude that the ultra-weak closure of $B$ is contained in $M_n(A_{\overline{\infty}})$.
\end{proof}

With respect to the original question we are now left with two questions which seem to be quite unrelated to each other.
\begin{enumerate}
\item Let $M$ be a von Neumann algebra and let $A$ be an ultra-weakly dense $*$-sub-algebra of $M$. Under what circumstances is $M^A=M$?
\item Let $\Gamma$ be a discrete group and let $\Delta: \Cz \Gamma \to \cU^e$ be a derivation. Under what circumstances is $\Delta$ closable?
\end{enumerate}

As the proofs show, the study of the algebra $M^A$ is not circumvented when working with continuous cohomology, i.e. continuous derivations. Although
there is no question about uniqueness of extensions, but extending to more than $M^A$ seems to be difficult.

\section{Continuous cohomology}

Let us give an ad hoc definition of the first continuous $L^2$-Betti number. It is easy to see that this continuous version of $L^2$-Betti number is
indeed the dimension of a continuous version of the first cohomology group. Denote by $\Der_c(A,\cU^e)$ the right $\cU^e$-module of
closable derivations.

\begin{definition}
Let $A \subset M$ be a $*$-subalgebra of a finite von Neumann algebra $M$. We define
\[\eta_{1}(A,\tau) = \dim_{\cU^e} \Der_c(A,\cU^e) + \beta_0^{(2)}(M) - 1.\]
\end{definition}

It is clear that always $0 \leq \eta_1(A,\tau) \leq \beta_1^{(2)}(A,\tau)$. The following theorem subsumes the results we have obtained in terms
of the quantity $\eta_1$.

\begin{theorem}
Let $A \subset M$ be a ultra-weakly dense $*$-subalgebra of a finite von Neumann algebra $M$. The following holds:
\[  \eta_1(M) \leq \eta_1(A),  \quad \mbox{and} \quad \eta_{1}(A) = \eta_1(M^A) . \]
\end{theorem}
The study of $\eta_1(M,\tau)$ is much easier than the study of $\beta_1^{(2)}(M,\tau)$, however, there is no example where we can show
that $\eta_1$ does not vanish.

\begin{theorem} Let $M$ and $N$ be finite von Neumann algebras. The following inequality holds:
\[ \eta_1(M \ast N) \leq \eta_1(M) - \beta_0^{(2)}(M)+ \eta_1(N)  - \beta_0^{(2)}(N) + 1.\]
\end{theorem}
\begin{proof}
Closable derivations are continuous by Theorem \ref{closcont}. We get that $\Der_c(M \ast N,\cU^e)$ can be identified with a sub-module of
$\Der_c(M,\cU^e) \oplus \Der_c(N,\cU^e)$ by restriction. The rest of the argument is straighforward.
\end{proof}

\begin{theorem} \label{top}
Let $M$ be a finite von Neumann algebra that either:
\begin{enumerate}
\item is a tensor product of two diffuse algebras,
\item admits a Cartan sub-algebra, or
\item has diffuse center.
\end{enumerate}
If any of the above conditions is satisfied, then $\eta_1(M)=0$.
\end{theorem}

We first proof the following lemma.
\begin{lemma}
Let $\Delta: M \to \cU^e$ be a derivation that vanishes on a diffuse sub-algebra $A \subset M$ and let $u \in M$ be a unitary operator which normalizes $A$. Under
these circumstances $\Delta(u) =0$.
\end{lemma}
\begin{proof}
Note that $\Delta(u^*) = - (u^* \otimes u^*) \Delta(u)$.
Let $h$ be a diffuse element in $A$.
\begin{eqnarray*} 0 &=& \Delta(uhu^*) \\
&=& (1 \otimes hu^*)\Delta(u) + (uh \otimes 1) \Delta(u^*)\\
&=& (1 \otimes hu^* - uhu^* \otimes u^*) \Delta(u) \\
&=& (u \otimes 1) (1 \otimes h - h \otimes 1) (u^* \otimes u^*) \Delta(u)
\end{eqnarray*}
Since $h$ is diffuse, $1 \otimes h -h \otimes 1$ is a non-zero divisor, hence $\Delta(u) =0$. This finishes the proof.
\end{proof}

\begin{proof}[Proof of Theorem \ref{top}]
$(1)$ Let $\Delta: M \to \cU^e$ be a continuous derivation and assume that $M$ decomposes as a tensor product $M = M_1 \hot M_2$. Pick
a diffuse abelian sub-algebra $A \subset M_1$. The derivation $\Delta$ is inner on $A \otimes 1$ and substracting the inner derivation
we get a derivation $\Delta'$ that vanishes on $A \otimes 1$. By the preceding lemma, $\Delta'$ has to vanish on $1 \otimes M_2$ and in
turn also on $M_1 \otimes 1$. By continuity we get that $\Delta'$ vanishes on $M$ and hence $\Delta$ was inner.

A similar argument applies in the case $(2)$ of existence of a Cartan sub-algebra $A \subset M$. The derivation minus some inner derivation vanishes 
on $A$ and hence on each unitary normalizing $A$ by the preceding lemma. By definition of Cartan subalgebra, these unitaries generate $M$. 
Whence the result follows.

The case $(3)$ of diffuse centre is a trivial consequence of the preceding lemma.
\end{proof}
\section{Conclusion}

The motivating question about the relation between the $L^2$-Betti numbers of a discrete group and its associated von Neumann algebra remains open. However,
we hope that some of the presented results will help to find a better understanding of this problem. 

It would be interresting to study the question about closability in the case of the derivations which arise as free difference quotients. In the case of
finite free Fisher information, it is known that these derivations are closable from $L^2(M,\tau)$ to $L^2(M \hot M^o, \tau \otimes \tau)$, see \cite{Voi5}. 
However, the domain of the adjoint of $\Delta$ is to small to derive the desired result directly.

\end{document}